\documentclass[a4paper,10pt]{article}
\usepackage{graphicx}
\usepackage{comment}
\usepackage{amssymb}
\usepackage{amsmath}
\usepackage{nicefrac}
\usepackage{graphicx}
\usepackage{dcolumn}
\usepackage{bm}
\usepackage{comment}
\usepackage{multirow}
\usepackage{cite}
\usepackage{xcolor}
\usepackage{url}
\usepackage{mathrsfs}
\usepackage{supertabular}
\usepackage[normalem]{ulem}
\usepackage{lscape}
\usepackage{soul}
\usepackage{subcaption} 
\usepackage{bold-extra}
\usepackage{hyperref}
\usepackage{latexsym,amsthm} 

\newtheorem{theorem}{Theorem}

\theoremstyle{definition}

\begin{document}
            
\huge

\begin{center}
A sum rule for derangements
\end{center}

\vspace{0.5cm}

\large

\begin{center}
Jean-Christophe Pain$^{a,b,}$\footnote{jean-christophe.pain@cea.fr}
\end{center}

\normalsize

\begin{center}
\it $^a$CEA, DAM, DIF, F-91297 Arpajon, France\\
\it $^b$Universit\'e Paris-Saclay, CEA, Laboratoire Mati\`ere en Conditions Extr\^emes,\\
\it 91680 Bruy\`eres-le-Ch\^atel, France\\
\end{center}

\vspace{0.5cm}

\begin{abstract}
We propose a sum rule for derangements. Three different proofs are provided. The first one involves integral representations and the second one relies on the Hermite identity for the integer part of the product of an integer by a real number. The third one, which, unlike the previous ones, is not based on induction, proceeds directly from a recurrence relation for the number of derangements.
\end{abstract}

\section{Introduction}\label{sec0}

A \textit{derangement} of $\mathscr{E}_n=\{1,2,3,\cdots,n\}$ is a permutation of $\mathscr{E}_n$ that has no fixed points. The number of derangements (sequence A000166 in OEIS \cite{sloane2023}) is given by
\begin{equation*}
    D(n)=n!\sum_{i=0}^{n}\frac{(-1)^i}{i!}.
\end{equation*}
One also has, for $n\geq 1$:
\begin{equation*}
    D(n)=\Big\lvert\Big\lvert\frac{n!}{e}\Big\lvert\Big\lvert=\Big\lfloor\frac{n!}{e}+\frac{1}{2}\Big\rfloor,
\end{equation*}
where $\lfloor x\rfloor$ represents the integer part of $x$ and $||x||$ represents the nearest integer to $x$, as well as
\begin{equation}\label{bas}
    D(n)=\left\lfloor {\frac {n!+1}{e}}\right\rfloor.
\end{equation}
Derangements have many interesting properties (see, for instance, \cite{apostol1974,graham1994,hassani2003,hassani2004}). Several identities involving derangements are known, such as, for an integer $p\geq 1$ (see, for instance, \cite{vinh2005})
$$
\sum_{n=0}^p\,\binom{p}{n}\,D(n)=p!,
$$
and, if $p$ and $l$ are integers ($p\geq 1$ and $l<p$)
$$
\sum_{n=0}^p\,\binom{p-l}{n-l}\,D(p-n)=(p-l)!.
$$
In the following, we provide an expression for the partial sum $\sum_{n=0}^p\,n\,D(n)$.

\begin{theorem}

Given an integer $p\geq 0$, we have

$$
\sum_{n=0}^p\,n\,D(n)=\Big\lfloor\frac{(p+1)!}{e}\Big\rfloor.
$$

\end{theorem}

Three proofs of Theorem 1 will be provided. The first one, given in Section \ref{sec1}, relies on integral representations; the second one, in Section \ref{sec2}, uses the Hermite identity, and the third one, detailed in Section \ref{sec3}, is a direct (non-inductive) proof involving  a recurrence relation. 

\section{First Proof Using Integral Representations}\label{sec1}

The proof of Theorem 1 given in this section is based on induction and involves integral representations.\\

{\em Proof of Theorem 1}\\
\\
Let us define
\begin{equation}\label{sp}
    S_p=\sum_{n=0}^p\,n\,D(n).
\end{equation}
We will prove the result 
$$
S_p=\Big\lfloor\frac{(p+1)!}{e}\Big\rfloor
$$
by induction. We obviously have
$$
S_0=0\;\;\;\;\mathrm{and}\;\;\;\;\Big\lfloor\frac{1!}{e}\Big\rfloor=0
$$
as well as
$$
S_1=D(1)=0\;\;\;\;\mathrm{and}\;\;\;\;\Big\lfloor\frac{2!}{e}\Big\rfloor=0,
$$
and, since $D(2)=1$,
$$
S_2=D(1)+2D(2)=2\;\;\;\;\mathrm{and}\;\;\;\;\Big\lfloor\frac{3!}{e}\Big\rfloor=2.
$$
Identity \eqref{sp} is thus valid for the first three values of $p$. Assuming it is true for $p$, one has to show, that for $p+1$,
$$
S_{p+1}=\sum_{n=0}^{p+1}\,n\,D(n)=\Big\lfloor\frac{(p+2)!}{e}\Big\rfloor,
$$
which is equivalent to
$$
S_{p+1}=S_p+(p+1)\Big\lfloor\frac{(p+1)!+1}{e}\Big\rfloor.
$$
In other words, using the expression for the number of derangements in Equation \eqref{bas}, one needs to prove
\begin{equation}\label{toprove}
\Big\lfloor\frac{(p+2)!}{e}\Big\rfloor=\Big\lfloor\frac{(p+1)!}{e}\Big\rfloor+(p+1)\Big\lfloor\frac{(p+1)!+1}{e}\Big\rfloor.
\end{equation}

Hassani \cite[Corollary 3.2]{hassani2003} found that
\begin{equation}\label{has1}
    e\Big\lfloor\frac{n!+1}{e}\Big\rfloor=\int_{-1}^{\infty}\,e^{-t}\,t^n\,\mathrm{d}t,
\end{equation}
as well as
$$
\int_{-1}^0\,e^{-t}\,t^n\,\mathrm{d}t=\left\{
\begin{array}{llll}
    -e\left\{\frac{n!}{e}\right\} & \mathrm{if} & n &\mathrm{odd},\\
    e-e\left\{\frac{n!}{e}\right\} & \mathrm{if} & n &\mathrm{even},
\end{array}
\right.
$$
which can be recast into
\begin{equation}\label{has2}
    \int_{-1}^0\,e^{-t}\,t^n\,\mathrm{d}t=e\left(\frac{1+(-1)^n}{2}-\left\{\frac{n!}{e}\right\}\right),
\end{equation}
where $\left\{x\right\}$ denotes the fractional part of $x$. Then, using $\left\{x\right\}=x-\lfloor x\rfloor$ and Equation \eqref{has2}, one obtains
\begin{equation}\label{diffe}
    \Big\lfloor\frac{(p+1)!}{e}\Big\rfloor=\frac{1}{e}\int_{-1}^0\,e^{-t}\,t^{p+1}\,\mathrm{d}t+\frac{(p+1)!}{e}-\frac{[1+(-1)^{p+1}]}{2}
\end{equation}
and
$$
\Big\lfloor\frac{(p+2)!}{e}\Big\rfloor=\frac{1}{e}\int_{-1}^0\,e^{-t}\,t^{p+2}\,\mathrm{d}t+\frac{(p+2)!}{e}-\frac{[1+(-1)^{p+2}]}{2}.
$$
The difference of the two latter equations reads
\begin{equation*}
    \Big\lfloor\frac{(p+2)!}{e}\Big\rfloor-\Big\lfloor\frac{(p+1)!}{e}\Big\rfloor=\frac{1}{e}\int_{-1}^0\,e^{-t}\,t^{p+1}(t-1)\,\mathrm{d}t+(-1)^{p+1}+\frac{(p+2)!}{e}-\frac{(p+1)!}{e},
\end{equation*}
i.e., since
\begin{equation}\label{fact}
    n!=\int_0^{\infty}\,e^{-t}\,t^n\,\mathrm{d}t,
\end{equation}
one has
\begin{equation}\label{las}
   \Big\lfloor\frac{(p+2)!}{e}\Big\rfloor-\Big\lfloor\frac{(p+1)!}{e}\Big\rfloor=\frac{1}{e}\int_{-1}^0\,e^{-t}\,t^{p+1}(t-1)\,\mathrm{d}t+(-1)^{p+1}+\frac{1}{e}\int_0^{\infty}\,e^{-t}\,t^{p+1}\,(t-1)\,\mathrm{d}t.
\end{equation}
Integrating by parts, one easily gets that the right-hand side of Equation \eqref{las} equals
\begin{equation}\label{iint}
    \int_{-1}^{\infty}\,e^{-t}\,t^{p+1}\,\left(t-(p+2)\right)\,\mathrm{d}t=(-1)^p\,e
\end{equation}
and thus, after multiplication by $e$ and expressing $(-1)^{p+1}$ using Equation \eqref{iint}, we obtain
\begin{align*}
     e\left(\Big\lfloor\frac{(p+2)!}{e}\Big\rfloor-\Big\lfloor\frac{(p+1)!}{e}\Big\rfloor\right)&=\int_{-1}^0\,e^{-t}\,t^{p+1}(t-1)\,\mathrm{d}t\nonumber\\
     &\quad-\int_{-1}^{\infty}\,e^{-t}\,t^{p+1}\,\left(t-(p+2)\right)\,\mathrm{d}t\nonumber\\
     &\quad+\int_0^{\infty}\,e^{-t}\,t^{p+1}(t-1)\,\mathrm{d}t\nonumber\\
    &=(p+1)\int_{-1}^{\infty}\,e^{-t}\,t^{p+1}\,\mathrm{d}t.
\end{align*}
Finally, using the representation of Equation \eqref{has1}, Equation \eqref{las} becomes:
$$
\Big\lfloor\frac{(p+2)!}{e}\Big\rfloor-\Big\lfloor\frac{(p+1)!}{e}\Big\rfloor=\frac{(p+1)}{e}\int_{-1}^{\infty}\,e^{-t}\,t^{p+1}\,\mathrm{d}t=(p+1)\Big\lfloor\frac{(p+1)!+1}{e}\Big\rfloor,
$$
which completes the proof. $\hfill\Box$

\vspace{5mm}

\section{Second Proof Using the Hermite Identity}\label{sec2}

As shown in Section \ref{sec1}, proving Theorem 1 amounts to proving Equation \eqref{toprove}. The proof presented in this section is inductive (like the previous one), but relies on the Hermite identity.\\

{\em Proof of Theorem 1}\\
\\
The Hermite identity \cite{matsuoka1964} reads
$$
\sum _{k=0}^{m-1}\left\lfloor x+{\frac {k}{m}}\right\rfloor =\lfloor mx\rfloor.
$$
This can be easily shown by noticing that the function
$$
f(x)=\lfloor x\rfloor +\left\lfloor x+{\frac {1}{m}}\right\rfloor +\ldots +\left\lfloor x+{\frac {m-1}{m}}\right\rfloor -\lfloor mx\rfloor
$$
is $1/m$-periodic. Indeed, since $f(x)=0$ for all $x\in [0,1/m[$, it follows that $f(x) = 0$ for all real $x$.

Let us set $m=p+2$ and $x=(p+1)!/e$. The Hermite identity becomes
\begin{align*}
\Big\lfloor\frac{(p+2)!}{e}\Big\rfloor&=\Big\lfloor\frac{(p+1)!}{e}\Big\rfloor+\Big\lfloor\frac{(p+1)!}{e}+\frac{1}{p+2}\Big\rfloor+\Big\lfloor\frac{(p+1)!}{e}+\frac{2}{p+2}\Big\rfloor\nonumber\\
&\quad+\cdots +\Big\lfloor\frac{(p+1)!}{e}+\frac{p+1}{p+2}\Big\rfloor.
\end{align*}
For $t\in[0,1]$, one has $-e<e^{-t}<e$ and thus
$$
-\frac{e}{p+2}<-e\frac{(-1)^p}{p+2}=-e\left[\frac{t^{p+2}}{p+2}\right]_0^{-1}=-e\int_0^{-1}t^{p+1}\,\mathrm{d}t<\int_0^{-1}t^{p+1}\,e^{-t}\,\mathrm{d}t,
$$
as well as
$$
\int_0^{-1}t^{p+1}\,e^{-t}\,\mathrm{d}t<e\int_0^{-1}t^{p+1}\,\mathrm{d}t=e\left[\frac{t^{p+2}}{p+2}\right]_0^{-1}=e\frac{(-1)^p}{p+2}<\frac{e}{p+2}
$$
and
$$
-\frac{e}{p+2}<\int_0^{\infty}t^{p+1}\,e^{-t}\,\mathrm{d}t-\int_{-1}^{\infty}t^{p+1}\,e^{-t}\,\mathrm{d}t<\frac{e}{p+2}.
$$
Using Equations \eqref{has1} and \eqref{fact}, one gets
$$
-\frac{1}{p+2}<\frac{(p+1)!}{e}-\Big\lfloor\frac{(p+1)!+1}{e}\Big\rfloor<\frac{1}{p+2},
$$
or
$$
-\frac{k}{p+2}<\frac{(p+1)!}{e}-\Big\lfloor\frac{(p+1)!+1}{e}\Big\rfloor<1-\frac{p+1}{p+2}\leq 1-\frac{k}{p+2}
$$
for $0\leq k\leq p$. This implies
\begin{equation}\label{solu}
    \Big\lfloor\frac{(p+1)!+1}{e}\Big\rfloor-\frac{(p+1)!+1}{e}<\frac{k}{p+2}-\frac{1}{e}<\Big\lfloor\frac{(p+1)!+1}{e}\Big\rfloor+1-\frac{(p+1)!+1}{e}.
\end{equation}
Let us now set $y=[(p+1)!+1]/e$ and $\alpha=k/(p+2)-1/e$. Since Equation \eqref{solu} means that $\lfloor y\rfloor-y<\alpha<\lfloor y\rfloor+1-y$, one has $\lfloor y+\alpha\rfloor=\lfloor y\rfloor$, i.e.,
$$
\Big\lfloor\frac{(p+1)!}{e}+\frac{k}{p+2}\Big\rfloor=\Big\lfloor\frac{(p+1)!+1}{e}\Big\rfloor
$$
and 
$$
\sum_{k=0}^{p}\Big\lfloor\frac{(p+1)!}{e}+\frac{k}{p+2}\Big\rfloor=(p+1)\Big\lfloor\frac{(p+1)!+1}{e}\Big\rfloor.
$$
This finally yields
$$
\Big\lfloor\frac{(p+2)!}{e}\Big\rfloor=\Big\lfloor\frac{(p+1)!}{e}\Big\rfloor+(p+1)\Big\lfloor\frac{(p+1)!+1}{e}\Big\rfloor,
$$
which completes the proof. $\hfill\Box$

\vspace{5mm}

\section{Direct (Non-Inductive) Proof Using a Recurrence Relation}\label{sec3}

In this section we present a non-inductive proof of Theorem 1 (i.e., which is not based on a proof of the identity of Equation \eqref{toprove}). The proof involves a recurrence relation for the number of derangements.\\

{\em Proof of Theorem 1}\\
\\
The derangements satisfy the well-known relation
\begin{equation}\label{recc1}
    D(n)=(n-1)(D(n-1)+D(n-2)).
\end{equation}
The relation \eqref{recc1} can be recast into
\begin{equation*}
    D(n)=nD(n-1)-(D(n-1)-(n-1)D(n-2)),
\end{equation*}
which gives
\begin{align*}
    D(n)-nD(n-1)&=-(D(n-1)-(n-1)D(n-2))\\
                &=\cdots\\
                &=(-1)^{n-2}(D(2)-2D(1))\\
                &=(-1)^n,
\end{align*}
the last equality stemming from the fact that $D(2)-2D(1)=1$, as seen in Section \ref{sec1}. We thus have 
\begin{equation*}
    D(n)-D(n-1)=(n-1)D(n-1)+(-1)^n,
\end{equation*}
for which different proofs and combinatorial interpretations were provided. Thus, one has
\begin{equation*}
    \sum_{n=2}^{p+1}(D(n)-D(n-1))=\sum_{n=2}^{p+1}(n-1)D(n-1)+\sum_{n=2}^{p+1}(-1)^n,
\end{equation*}
or
\begin{equation*}
    D(p+1)=\sum_{n=0}^pnD(n)+\frac{1-(-1)^{p}}{2},
\end{equation*}
yielding
\begin{equation*}
    S_p=D(p+1)-\frac{1-(-1)^{p}}{2},
\end{equation*}
or also
\begin{equation*}
    S_p=\Big\lvert\Big\lvert\frac{(p+1)!}{e}\Big\lvert\Big\lvert-\frac{1-(-1)^{p}}{2}.
\end{equation*}
Using Equation \eqref{diffe}, one gets
\begin{equation*}
    \frac{(p+1)!}{e}-\Big\lfloor\frac{(p+1)!}{e}\Big\rfloor=-\frac{1}{e}\int_{-1}^0\,e^{-t}\,t^{p+1}\,\mathrm{d}t+\frac{[1+(-1)^{p+1}]}{2}.
\end{equation*}
Setting 
$$
f(p)=-\frac{1}{e}\int_{-1}^0\,e^{-t}\,t^{p+1}\,\mathrm{d}t=\frac{(-1)^p}{e}\int_{0}^1\,e^{t}\,t^{p+1}\,\mathrm{d}t,
$$
one sees that, for $p$ even, $f(p)$ is a positive decreasing function of $p$ and since $f(0)=1/e<1/2$, we obtain that $f(p)<1/2$. On the other hand, when $p$ is odd, then  $f(p)$ is a negative increasing function of $p$ and since $f(1)=(2-e)/e$, it follows that $f(1)+1=2/e>1/2$ and $f(p)>1/2$. Thus,
$$
-\frac{1}{e}\int_{-1}^0\,e^{-t}\,t^{p+1}\,\mathrm{d}t+\frac{[1+(-1)^{p+1}]}{2}
\begin{cases}
    <1/2 & \text{if } p \text{ is even},\\\\
    >1/2 & \text{if } p \text{ is odd}, 
\end{cases}
$$
and then
$$
\frac{(p+1)!}{e}
\begin{cases}
    <\displaystyle\Big\lfloor\frac{(p+1)!}{e}\Big\rfloor+1/2 & \text{if } p \text{ is even},\\\\
    >\displaystyle\Big\lfloor\frac{(p+1)!}{e}\Big\rfloor+1/2 & \text{if } p \text{ is odd}.
\end{cases}
$$
This gives
\begin{equation*}
\Big\lvert\Big\lvert\frac{(p+1)!}{e}\Big\lvert\Big\lvert=
\begin{cases}
  \displaystyle\Big\lfloor\frac{(p+1)!}{e}\Big\rfloor & \text{if } p+1 \text{ is odd},\\\\
  \displaystyle\Big\lfloor\frac{(p+1)!}{e}\Big\rfloor + 1 & \text{if } p+1 \text{ is even},
\end{cases}
\end{equation*}
yielding finally 
\begin{equation*}
    S_p=\Big\lfloor\frac{(p+1)!}{e}\Big\rfloor,
\end{equation*}
which completes the proof. $\hfill\Box$

\vspace{5mm}

It is worth mentioning that in an article about permutation enumeration \cite{sedgewick1977}, Sedgewick proposed an improvement of Heap's algorithm for the generation of permutations by interchanges \cite{heap1963}. It turns out that the complexity of Sedgewick's efficient algorithm involves the quantity $A_N$ defined by the induction relation  
\begin{equation*}
    A_N=N A_{N-1}+\begin{cases}
    0 & \text{if } N \text{ is even},\\
    N-1 & \text{if } N \text{ is odd},\\
    \end{cases}
\end{equation*}
for $N>1$ and $A_1=0$. One actually has
\begin{equation*}
    A_N=N!~\sum_{k=2}^N\frac{(-1)^k}{k!}=\left\lfloor\frac{N!}{e}\right\rfloor,
\end{equation*}
and there is a strong connection between the sum $S_p$ and the quantity $A_N$, which are related to each other by
\begin{equation*}
    S_p=A_{p+1}.
\end{equation*}

\section{Conclusion}

We obtained a sum rule for derangements, proven by three different methods. The first two methods are inductive proofs; the first one uses integral representations of the integer part of $(n!+1)/e$ and of the fractional part of $n!/e$ ($n$ being an integer), and the second one is based on the Hermite identity, expanding the integer part of $nx$ (where $x$ is a real number) as a sum of integer parts of $x+k/n$, with $k$ ranging from 0 to $n-1$. The third proof proceeds directly from a recurrence relation for derangements. The identity presented in Theorem 1 may be useful for deriving sum rules for Laguerre polynomials \cite{even1976}.


\end{document}